\documentclass{amsart}
\usepackage{amsthm}
\usepackage{amsmath}
\usepackage{amstext}
\usepackage{amsfonts}
\usepackage{latexsym}
\usepackage{amscd}
\usepackage{xypic}
\newcommand{\la}{\ensuremath{\longrightarrow}}

\newcommand{\sheaf}{\ensuremath{\mathcal{O}}}

\theoremstyle{plain}
\newtheorem{theorem}{Theorem}[section]
\newtheorem{lemma}[theorem]{Lemma}
\newtheorem{proposition}[theorem]{Proposition}

\theoremstyle{definition}
\newtheorem{definition}[theorem]{Definition}

\newtheorem{example}[theorem]{Example}
\newtheorem{remark}[theorem]{Remark}
\newtheorem{question}[theorem]{Question}
\numberwithin{equation}{section}

\begin{document}

\title[Divisorial extremal neighborhoods over $cE_7$ and $cE_6$ singularities]{3-fold divisorial extremal neighborhoods over  $cE_7$
and $cE_6$ Compound DuVal Singularities}
\author{Nikolaos Tziolas}
\address{Department of Mathematics, University of Cyprus, P.O. Box 20537, Nicosia, 1678, Cyprus}
\email{tziolas@ucy.ac.cy}




\begin{abstract}
Let $P \in \Gamma \subset X$ be the germ of a Gorenstein 3-fold singularity and $\Gamma$ a smooth curve through it such that 
the general hyperplane section $S$ of $X$ containing $\Gamma$ is DuVal of type $E_6$ or $E_7$. In this paper we obtain criteria 
for the existence of a terminal divisorial extremal neighborhood $C \subset Y \stackrel{f}{\la} X \ni P$ contracting an irreducible 
divisor $E$ onto $\Gamma$ and classify all such neighborhoods. 
\end{abstract}

\maketitle

\section{Introduction}

The study of the structure of birational maps between two Mori fiber spaces is one of the most interesting open problems in 3-fold birational geometry. 
In recent years a lot of progress has been made in this direction especially since the introduction of the Sarkisov program by Corti, Reid and 
Sarkisov~\cite{Cor95}~\cite{Cor-Rei00}. 
The aim of this program is to factor any birational map between two Mori fiber spaces as the composition of simpler maps called elementary links. 
Each link is composed of flips, inverse flips 
and divisorial contractions. In this context it is important to understand such birational maps. The structure of divisorial contractions is still 
unknown in its complete generality, even though 
lately a lot of progress has been made in this direction by Corti, Cutkosky, Kawamata, Kawakita, Mori,Tziolas and others. 

A divisorial contraction is a proper birational map $f \colon E \subset Y \la X \supset \Gamma$ such that $X$ and $Y$ have at worst terminal 
singularities, $-K_Y$ is $f$-ample, 
$E$ is an irreducible divisor, $Y-E\cong X-\Gamma$ and $\Gamma=f(E)$. Given $\Gamma \subset X$, one would like to know when there exists a 
divisorial contraction contracting an irreducible 
divisor $E$ onto $\Gamma$. There are two main cases to consider here. The first one is when $\Gamma$ is just a point. In this direction there 
has been many results by Corti~\cite{Cor-Rei00}, 
Kawakita~\cite{Ka01},~\cite{Ka02} and Kawamata~\cite{Kaw94}. The second case is when $\Gamma$ is a curve. In this direction there has been 
results by Kawamata~\cite{Kaw94} and 
Tziolas~\cite{Tzi03}~\cite{Tzi05a}~\cite{Tzi05}. Kawamata has shown that if there is a cyclic quotient singularity of $X$ on $\Gamma$ then 
there is no divisorial contraction contracting an 
irreducible divisor onto $\Gamma$. 

Given a curve $\Gamma$ in a 3-fold $X$, the existence of a divisorial contraction  $f \colon E \subset Y \la X \supset \Gamma$ is a 
local question around the singularities of $X$ on $\Gamma$~\cite[Proposition 1.2]{Tzi03}. 
Therefore the proper setting of the problem is to replace $X$ by an analytic neighborhood of a singular point $P\in X$. 
In this context, a divisorial 
contraction is the same as a divisorial extremal neighborhood as introduced by Koll\'ar and Mori~\cite{Ko-Mo92}.

By Reid's general elephant conjecture, given the germ $P \in \Gamma \subset X$ of a 3-fold terminal singularity and a curve through it, 
if a contraction exists then the 
general member of $|-K_X|$ containing $\Gamma$ 
has DuVal singularities~\cite{Ko-Mo92}. This can be of type $A_n$, $D_n$, $E_6$, $E_7$ and $E_8$. In the series of 
papers~\cite{Tzi03}~\cite{Tzi05a}~\cite{Tzi05} 
a complete classification was given in the cases when $X$ is Gorenstein, $\Gamma$ a smooth curve and the general hyperplane 
section of $X$ through $\Gamma$ is DuVal of type $A_n$ and $D_n$. 
For the classification to be complete in the Gorenstein case, it remains to treat the exceptional cases. The $E_8$ does not 
appear since $E_8$ is factorial and $\Gamma$ is assumed 
to be smooth. 

In this paper we present a complete classification of the cases when the general hyperplane section of $X$ through $\Gamma$ is DuVal of 
type $E_6$ or $E_7$. The method is mostly a combination 
of the methods that were used in~\cite{Tzi03},~\cite{Tzi05a} and~\cite{Tzi05}. So, as far as methodology is concerned there is little new in this paper. 
However, I hope that the results 
will be of interest to the people that work in the explicit 3-fold Minimal Model Program and especially the Sarkisov program.

Let $P \in X$ be the germ of a 3-fold singularity and $\Gamma$ a smooth curve through it.
Theorems~\ref{E7},~\ref{E6} give a complete classification of terminal divisorial extremal neighborhoods 
$C \subset Y \stackrel{f}{\la} X \ni P$ in the cases when the general hyperplane section $S$ of $X$ through $\Gamma$ 
is DuVal of type $E_7$ and $E_6$ respectively. In particular, it is shown that if the general hyperplane section $H$ of $X$ 
is DuVal of type $E_7$ or $D_n$ with $n\geq 6$, then there is no such neighborhood. In the remaining cases when $H$ is DuVal of 
type $E_6$, $D_5$ or $D_4$, explicit conditions for the 
existence of a neighborhood are given in terms of the equations of $P \in \Gamma \subset X$. 

The classification is done by explicit calculations using normal forms for the equation of $P\in X$ and the construction 
that first appeared in~\cite{Tzi03}. Normal forms are obtained in Proposition~\ref{normal-forms}.

Proposition~\ref{deformations} shows that there is a one to one correspondence between divisorial extremal neighborhoods 
$f\colon C \subset Y \la X \ni P$ and $\mathbb{Q}$-Gorenstein 
terminal smoothings of certain surface germs $C \subset Z$. 

Finally I must mention that some of the results that appear in this paper can be found in~\cite[Section 4]{Ko-Mo92}. 
In that paper, by using different methods, a classification of index 2 extremal neighborhoods $f\colon C \subset Y \la X \ni P$ 
is given. The classification of these neighborhoods is done by describing all the possible $Y$'s. This approach is not very useful 
from the point of view 
of the Sarkisov program since in this context, given a germ $\Gamma\subset X$, one would like to know when a divisorial 
contraction $f \colon E\subset Y \la X\supset \Gamma$ exists and 
to this question that method gives no answer. Moreover, given a germ $\Gamma\subset X$ as before such that the general 
hyperplane section of $X$ through $\Gamma$ is DuVal 
of type $E_7$, then if the contraction exists it is not known apriori that $Y$ has index 2 and so we cannot use that method 
directly to describe the neighborhoods.

\section{Construction of contractions.}
As was mentioned in the introduction, the existence of a terminal contraction is a local question around the singularity of $X$. 
The local version of a divisorial contraction
is that of an extremal neighborhood.

\begin{definition}[~\cite{Ko-Mo92}]
A three-dimensional extremal neighborhood is a proper birational morphism $C \subset Y \stackrel{f}{\la} X \ni P $, with the following properties:
\begin{enumerate}
\item $C \subset Y$ is the germ of a normal complex space along a proper curve $C$, $P \in X$ is the germ of a normal 3-fold singularity 
and $f^{-1}(P)_{red}=C$.
\item $-K_Y$ is $\mathbb{Q}$-Cartier and $f$-ample.
\end{enumerate}
\end{definition}
If $P\in X$ and $Y$ are both terminal, the neighborhood is called terminal.
If the exceptional set of $f$ is an irreducible divisor, then the extremal neighborhood is called divisorial. Otherwise, i.e.,
if it is one-dimensional, it is called flipping.

Let $P\in \Gamma \subset X$ be the germ of a Gorenstein terminal 3-fold singularity along a smooth curve $\Gamma$ such that the general hyperplane
section $S$ of $X$ through $\Gamma$ is DuVal.
Then according to~\cite{Tzi03}, there exists a divisorial extremal neighborhood $C \subset Y \stackrel{f}{\la} X \ni P$ such that $Y$ has
canonical singularities and the exceptional set of $f$ is an irreducible divisor $E$ that is contracted onto $\Gamma$.
Moreover, the contraction
is unique and it fits in the following commutative diagram
\begin{equation}\label{diagram}
\xymatrix{
Z \ar[d]_{h} \ar@{-->}[rr]^{\phi} & & Z^{\prime} \ar[d]^{\pi}\\
W \ar[dr]_g &   & Y \ar[dl]^f \\
            & X & }
\end{equation}
which is obtained as follows. $W$ is the blow up of $X$ along $\Gamma$. There are two $g$-exceptional divisors. A ruled surface $E$
over $\Gamma$ and $F\cong \mathbb{P}^2$ over $P$.
$Z$ is the $\mathbb{Q}$-factorialization of $E$ and $\phi$ a composition of flips. Finally $\pi$ is the contraction of the birational
transform $F^{\prime}$ of $F$ in $Z^{\prime}$.

In order to classify the divisorial extremal neighborhoods of interest in this paper, 
we will obtain normal forms for the equations of $P \in \Gamma \subset X$ and we 
will calculate explicitly all steps of the
above construction. In particular we will show that
the birational transform $F_Z$ of $F$ in $Z$ is still $\mathbb{P}^2$ and hence there are no flips and $\pi $ contracts $F_Z$ to a
terminal singularity. Then $Y$ is terminal if and only if the singularities of $Z$ on the $h$-exceptional curves and not on $F_Z$, are terminal.

The most difficult part of the previous construction is the description of $Z$.
The $\mathbb{Q}$-factorialization of a divisor in a terminal 3-fold is known to exist but in general its explicit description is hard.
The easiest case is when $W$ has isolated hypersurface singularities on $E$, and then $Z$ is just the blow up of $W$ along $E$.
Otherwise $Z$ is the blow up of a multiple of $E$. To determine which multiple has to be blown up, it is necessary to study the
singularities of $W$ at the generic points of $E\cap F$. This is not hard and can be done as follows.
Let $S \in |-K_X|$ be the general hyperplane section of $X$ through $\Gamma$. Then $S=g_{\ast}S_W$, where $S_W\in |-K_W|$ is the
general section, and $g \colon S_W \la S$ is
the blow up of $S$ along $\Gamma$. Therefore, the extended dual graph of $E \subset W$ at the generic points of $E\cap F$ is the
same as the extended dual graph of $\Gamma_W \subset S_W$,
where $\Gamma_W$ is the birational transform of $\Gamma$ in $S_W$. In the cases of interest of this paper, $S$ is DuVal of type either $E_7$ or $E_6$.

If $S$ is of type $E_7$, then in suitable local analytic coordinates $S$ is given by $x^2+y^3+yz^3=0$ and $\Gamma$ by $x=y=0$.
A straightforward calculation shows that
$S_W$ is smooth along $\Gamma_W$ and hence $W$ has isolated singularities along $E$. Therefore in this case $Z$ is just the blow up of $E$.

If $S$ is of type $E_6$, then in suitable local analytic coordinates $S$ is given by $x^2+y^3-z^4=0$ and $\Gamma$ by $x-z^2=y=0$.
A straightforward calculation shows that
$S_W$ is DuVal of type $D_5$ given by $x^2+y^2z-z^4=0$ and $\Gamma_W$ is given by $x=z=0$. Therefore $W$ is singular along $E\cap F$ and
moreover $2E$ is Cartier generically along $E\cap F$.
Hence $Z$ is the blow up of $2E$.

\section{Normal Forms.}
The next proposition obtains normal forms for the equations of $P\in \Gamma \subset X$.
\begin{proposition}\label{normal-forms}
Let $P \in \Gamma \subset X$ be the germ of a 3-fold Gorenstein singularity $0\in X$ along a smooth curve $\Gamma$. 
Let $S$ be the general hyperplane section of $X$ 
that contains $\Gamma$ and $H$ the general hyperplane section of $X$. Then
\begin{enumerate}
\item Suppose that $S$ is DuVal of type $E_7$. Then in suitable local analytic coordinates, the equation of $P \in X$ is
\[
x^2+y^3+yz^3+tf_2(y,t)+tf_{\geq 3} (y,z,t)=0
\]
and $\Gamma$ is given by $x=y=t=0$. $f_2(y,t)$ is a homogeneous polynomial of degree 2 in $y$ and $t$,  $f_{\geq 3}(y,z,t)$ 
has lowest degree at least 3 and $f_{\geq 3}(0,z,0)=0$.
Moreover, 
\begin{enumerate}
\item $H$ is DuVal of type $D_n$ if and only if the cubic term $y^3+tf_2(y,t)$ is not a cube. If $n \geq 6$ then  $f_2(0,t)=0$. 
Suppose that  $f_2(0,t)\neq 0$. Then
\begin{enumerate}
\item $H$ is of type $D_4$ if and only if the cubic term $g=y^3+tf_2(y,t)$ is square free, i.e., it is a product of three independent linear forms.
\item $H$ is of type $D_5$ if and only if the cubic term $g=y^3+tf_2(y,t)$ is divisible by a square, i.e., $g=l_1^2(y,t)l_2(y,t)$, 
where $l_1(y,t)$, $l_2(y,t)$ are independent linear
forms.
\end{enumerate}
\item $H$ is DuVal of type $E_6$ if and only if $y^3+tf_2(y,t)$ is a cube and either $f_2(y,t)\neq 0$ (equivalently, $f_2(0,t)\neq 0$) or $t^3$ 
appears in $f_{\geq 3}(y,z,t)$.
\item $H$ is DuVal of type $E_7$ if and only if $y^3+tf_2(y,t)$ is a cube, $f_2(y,t)= 0$ (equivalently, $f_2(0,t)= 0$) 
and there is no $t^3$ in $f_{\geq 3}(y,z,t)$.
\end{enumerate}
\item Suppose that $S$ is DuVal of type $E_6$. Then in suitable local analytic coordinates the equation of $P \in X$ is
\[
x^2+xz^2+y^3+ayt^2+by^2t+tf_{\geq 3}(y,z,t)=0
\]
where $\Gamma$ is given by $x=y=t=0$ and no monomials of the form $z^k$, $k\geq 4$, or $y^s$, $s \geq 3$, appear in $f_{\geq 3}(y,z,t)$.
Moreover, 
\begin{enumerate}
\item If $H$ is DuVal of type $D_4$, then $a\neq 0$ and $b^2-4a \neq 0$.
\item If $H$ is DuVal of type $D_5$, then $a=0$ and $b \neq 0$.
\item If $H$ is DuVal of type $E_6$, then $a=b=0$.
\end{enumerate}
\end{enumerate}
\end{proposition}
\begin{proof}
Suppose that $S$ is DuVal of type $E_7$. Then~\cite{Jaf92} in suitable local analytic coordinates $S$ is given by
$x^2+y^3+yz^3=0$ and $\Gamma$ by $x=y=0$. Hence $X$ is given by an equation of the form
\[
x^2+y^3+yz^3+tf(x,y,z,t)=0
\]
and $\Gamma$ by $x=y=t=0$. By using the Weierstrass preparation theorem it takes the form
\[
x^2+y^3+yz^3+tf(y,z,t)=0
\]
Moreover, if $f(y,z,t)$ has a linear term, then $X$ is of type $cA_n$, for some $n$. Hence the general hyperplane section of $X$ 
containing $\Gamma$ is of type $A_m$, for some $m$~\cite{Tzi03}.
Let $f_2(y,z,t)$ be the degree 2 homogeneous part of $f(y,z,t)$. Suppose $f_2(y,z,t)=a_1y^2+a_2t^2+a_3yt+a_4yz+a_5zt+a_6z^2$. Then

\textit{Claim:} If $(a_4,a_5,a_6)\neq (0,0,0)$, then the general hyperplane section of $X$ through $\Gamma$ is of type $D_n$.

To show the claim, assume that $(a_4,a_5,a_6)\neq (0,0,0)$ and let $S_{\lambda}$ be the section given by $y=\lambda t$. For general $\lambda$ it 
is normal because $t=0$ is normal.
Then $S_{\lambda}$ is given by
\[
x^2+(\lambda^3+a_1\lambda^2+a_3\lambda+a_2)t^3+(a_4\lambda+a_5)zt^2+a_6z^2t+\lambda tz^3 +tf_{\geq 3}(\lambda t,z,t)=0
\]
If  $(a_4,a_5,a_6)\neq (0,0,0)$, then the cubic term of the above equation is not a cube and therefore it is a DuVal singularity of type $D_n$, 
for some $n$, as shown by the next lemma.
\begin{lemma}[~\cite{AGZV85}]\label{Dn}
Let $0\in X$ be a normal surface singularity given by \[
x^2+f_{\geq 3}(y,t)=0
\]
Then $0\in X$ is DuVal of type $D_n$, for some $n$, if and only if the cubic term $f_3(y,t)$ of $f_{\geq 3}(y,t)$ is nonzero and not a cube.
\end{lemma}
Therefore the equation of $0\in X$ becomes
\begin{equation}
x^2+y^3+yz^3+tf_2(y,t)+tf_{\geq 3}(y,z,t)=0
\end{equation}
Suppose that $f_{\geq 3}(0,z,0)\neq 0$. Then the equation of $0\in X$ becomes
\[
x^2+y^3+yz^3+tf_2(y,t)+tf_{\geq 3}(y,z,t)+tz^nu=0
\]
where $n\geq 3$ and $u$ is a unit. The change of variables $y \mapsto y-tz^{n-3}u$ makes the equation of $0\in X$ as in $(1)$ with the 
added property that $f_{\geq 3}(0,z,0)=0$.

Next we study the general hyperplane section $H$ of $X$. It is DuVal of type $D_n$, $E_6$ or $E_7$. We will investigate when this 
happens according to the form of the equation of $X$.

\textit{Case 1:} Suppose that the cubic term of the equation of $X$, $y^3+tf_2(y,t)$ is not a cube.

Let $S_{\lambda}$ be the section
given by $z=\lambda t$. For general $\lambda$ it is normal because $t=0$ is normal. Then $S_{\lambda}$ is given by an equation of the form
\[
x^2+f_3(y,t)+f_{\geq 4}(y,t)
\]
where $f_3(y,t)$ is a homogeneous polynomial of degree 3 that is not a cube. Hence by Lemma~\ref{Dn},  $S_{\lambda}$ is DuVal of type $D_n$, for some $n$.

Now suppose that $f_2(0,t)\neq 0$. We will show that $n\leq 5$. There are two cases to consider.

The first case is when $y^3+tf_2(y,t)=l_1(y,t)l_2(y,t)l_3(y,t)$, where $l_i(y,t)$ are independent linear forms. In this case $H$ is DuVal of type $D_4$.

The only nontrivial case is when $y^3+tf_2(y,t)$ is not a cube and it is divisible by a square. In this case we will show that the 
general section $H$ of $X$ is DuVal of type $D_5$. The proof is based to the following elementary result.

\begin{lemma}
Let $P\in S$ be a DuVal surface singularity of type $D_n$, $E_6$ or $E_7$. Let $f \colon S^{\prime} \la S$ be the blow up of $P \in S$. 
Then the singularities of $S^{\prime}$ are
\begin{enumerate}
\item In the case that $P \in S$ is DuVal of type $D_n$,
\begin{enumerate}
\item One DuVal of type $A_1$ and one of type $D_{n-2}$, if $n \geq 6$
\item One DuVal of type $A_1$ and one of type $A_3$, if $n=5$
\item Three DuVal of type $A_1$, if $n=4$.
\end{enumerate}
\item One DuVal of type $A_5$ if $P\in S$ is of type $E_6$
\item One DuVal of type $D_6$ if $P\in S$ is of type $E_7$.
\end{enumerate}
\end{lemma}

Suppose now that $f_2(0,t)\neq 0$ and that the cubic term of the equation of $X$ is not a cube and is divisible by a square. Then
\[
 y^3+tf_2(y,t)=(a_1y+b_1t)^2(c_1y+d_1t)
 \]
with $a_1b_1c_1d_1\neq 0$. Then after a linear change of coordinates, the equation of $X$ takes the form
\[
x^2+y^2(ay+bt)+yz^3+c z^3t +tf_{\geq 3}(y,z,t)=0
\]
with $abc \neq 0$. Now let $H$ be the section given by $x=y+z+dt$. Its equation is
\[
(y+z+dt)^2+y^2(ay+bt)+yz^3+c z^3t +tf_{\geq 3}(y,z,t)=0
\]
I claim that it is DuVal of type $D_5$. This will follow from Lemma 3.3. So let $H^{\prime}$ be the blow up of $H$. In the affine chart given by $z=zt$, $y=yt$, $H^{\prime}$ is
given by
\[
(y+z+d)^2+y^2(ayt+bt)+yz^3t^2+cz^3t^2+\frac{1}{t}f_{\geq 3}(yt,zt,t)=0
\]
The change of variables $z \mapsto z-d$ transforms the above equation to
\[
(y+z)^2+y^2(ayt+bt)+y(z-d)^3t^2+c(z-d)^3t^2+\frac{1}{t}f_{\geq 3}(yt,(z-d)t,t)=0
\]
The degree 2 term is of the form $(y+z)^2+\beta t^2$, with $\beta \neq 0$. This is a product of two linear forms and therefore $H^{\prime}$ is DuVal of type $A_n$, with $n \geq 2$. Now 
according to lemma 3.3 the only case that this happens is when $H$ is $D_5$.

\textit{Case 2:} The cubic term $y^3+tf_2(y,t)$ is a cube.

Then the equations of $\Gamma \subset X$ become
\[
x^2+(y+bt)^3+yz^3+tf_{\geq 4}(y,z,t)=0
\]
By using the same method as in the previous case we will show that the general hyperplane section $H$ of $X$ given by $x=cy+dz+et$ is either $E_6$ or $E_7$. 
$H$ is given by the equation
\[
(cy+dz+et)^2+(y+bt)^3+yz^3+tf_{\geq 3}(y,z,t)=0
\]
Let $S^{\prime}\la S$ be the blow up of $S$. In the chart $y=yt$, $z=zt$, it is given by
\[
(cy+dz+e)^2+(y+b)^3t+yz^3t^2+\frac{1}{t}f_{\geq 3}(yt,zt,t)=0
\]
and the exceptional curve $C$ is $t=cy+dz+e=0$. We are looking for the singularities of $B_0H$ on $C$. A straighforward calculation shows that there is only one for
$t=cy+dz+e=y+b=0$. Setting $y \mapsto y-b$ and $z\mapsto z+z_0$, where $z_0=(-e+bc)/d$ the equation becomes
\begin{equation}
(cy+dz)^2+y^3t+(y-b)(z+z_0)^3t^2+\frac{1}{t}f_{\geq 3}(yt,zt,t)=0
\end{equation}
The degree 2 term of the equation is $(cy+dz)^2+\lambda t^2$. This is either a square or it factors. In any case it does not correspond to an $A_1$ point.
Similar calculations in the other charts show that there is no $A_1$ point and hence $H$ is not of type $D_n$, for any $n$. Therefore it is either $E_6$ or $E_7$.
If it is $E_7$ then $b=0$ because if $b\neq 0$, then for general choice of $c,\; d$ and $ e$, $\lambda \neq 0$ and hence $(cy+dz)^2+\lambda t^2$ factors. Hence $B_0H$ has an
$A_n$ singularity and hence it must be $E_6$.

Suppose now that $S$ is DuVal of type $E_6$.
Arguing in exactly the same way as in the $E_7$ case, we see that in suitable local analytic coordinates $X$ is given by
\[
x^2+y^3-z^4+tf_2(y,t)+t\phi (y,z,t)=0
\]
where $\Gamma$ is given by $x-z^2=y=t=0$ and $\mathrm{mult}_0(\phi ) \geq 3$.
Write $\phi(y,z,t)=z^ku_1+y^su_2+ f(y,z,t)$, with $u_1$ and $u_2$ units, $k\geq 4$ and $s \geq 3$. Then the equation of $X$ becomes
\[
x^2+y^3(1+ty^{s-3}u_1)-z^4(1-tz^{k-4}u_2)+tf_2(y,t)+tf(y,z,t)=0
\]
After the change of coordinates $y=y/\sqrt[3]{1+ty^{s-3}u_1}$ and $z=z/\sqrt[4]{1-tz^{k-4}u_2}$, the equation of $X$ takes the form
\[
x^2+y^3-z^4+tf_2(y,t)+tf(y,z,t)=0
\]
where $\mathrm{mult}_0(f)\geq 3$ and no monomial of the form $y^s$, $z^k$ with $s\geq 3$, $k\geq 4$, appear in $f(y,z,t)$.

Now suppose that the cubic term of the equation, $y^3+tf_2(y,t)$ is a cube. Then the general hyperplane section $H$ of $X$ is DuVal of type $E_6$.
Write $y^3+tf_2(y,t)=(y+at)^3$. then after a linear change of coordinates the equation
of $X$ takes the form
\[
x^2+y^3-z^4+tf(y,z,t)=0
\]
with $\mathrm{mult}_0(f)\geq 3$ and no $z^k$, $y^s$, $k\geq 4$, $s\geq 3$, appear in $f(y,z,t)$.

Suppose that $y^3+tf_2(y,t)$ is not a cube. Then if it is a product of three independent linear forms, $H$ is DuVal of type $D_4$ and after a 
suitable change of coordinates,
$X$ is given by \[
 x^2+y^3+ayt^2+by^2t-z^4+tf(y,z,t)=0
 \]
 with $a \neq 0$ and $b^2-4a\neq 0$.

 Finally suppose that $y^3+tf_2(y,t)$ is divisible by a square and not a cube. In this case we will show that $H$ is DuVal of type $D_5$. Write
 $y^3+tf_2(y,t)=(a_1y+b_1t)^2(c_1y+d_1t)$, with $a_1b_1c_1d_1 \neq 0$. After a linear change of coordinates the equation of $X$ becomes
 \[
 x^2+y^3+by^2t-z^4+tf(y,z,t)=0
 \]
 with $b \neq 0$. Suppose that $H$ is given by $x=y+z+\delta t$. Let $H^{\prime}$ be the blow up of $H$. Then arguing as in 
the $E_7$ case we see that $H^{\prime}$ is DuVal of type $A_r$, $r\geq 2$ and therefore $H$ must be $D_5$.

 Finally the change of variables $x=x-z^2$ brings the equation of $X$ in the claimed form.
\end{proof} 

\section{The $E_7$ Case.}

\begin{theorem}\label{E7}
Let $P\in \Gamma \subset X$ be the germ of a Gorenstein terminal 3-fold singularity along a smooth curve $\Gamma$ such that the general hyperplane
section $S$ of $X$ through $\Gamma$ is a DuVal
singularity of type $E_7$. Let $H$ be the general hyperplane section through $P$. Then $H$ is DuVal of type either $E_7$, $E_6$ or $D_n$, and
\begin{enumerate}
\item If $H$ is DuVal of type $E_7$ or $D_n$ with $n \geq 6$, then there is no terminal divisorial extremal neighborhood
$C \subset Y \stackrel{f}{\la} X \ni P$ from a
terminal 3-fold $Y$ to $X$, contracting
an irreducible divisor $E$ onto $\Gamma$.
\item In the remaining cases, when $H$ is either of type $E_6$, $D_5$ or $D_4$, in suitable local analytic coordinates 
$P\in X$ is given by\[
x^2+y^3+yz^3+tf_2(y,t)+tf_{\geq 3}(y,z,t) =0
\]
$\Gamma$ is given by $x=y=t=0$, 
$f_2(y,t)$ is a homogeneous polynomial of degree 2 in $y$ and $t$, and the lowest degree of $f_{\geq 3}(y,z,t)$ is at least 3.
Then a terminal divisorial extremal neighborhood
$C \subset Y \stackrel{f}{\la} X \ni P$ exists if and only if $f_2(0,t)\neq 0$.
In that case, let $S_Y \in |-K_Y|$ and $H_Y \in |\sheaf_Y |$ be the general members. Then
\begin{enumerate}
\item $C=f^{-1}(P)_{red}\cong \mathbb{P}^1$ and  $K_Y \cdot C = -1/2$.
\item $S_Y \cong S$.
\item $H_Y$ is normal.
\end{enumerate}
Moreover,
\begin{enumerate}
\item If $H$ is of type $E_6$, then the singular locus of $H_Y$ consists of exactly one
non log-canonical point of index 2 which is analytically isomorphic to
\[
(x^2+(y+t)^3y+y^2t^4=0)/\mathbb{Z}_2(1,1,1)
\]
and the extended dual graph of the minimal resolution of $C \subset H_Y$ is
\[
\xymatrix @R=5pt @C=5pt{
\underset{}{\overset{-2}{\circ}} \ar@{-}[r]   & \underset{}{\overset{-2}{\circ}} \ar@{-}[r] 
& \underset{}{\overset{-2}{\circ}} \ar@{-}[r] \ar@{-}[d] & \underset{}{\overset{-2}{\circ}} \ar@{-}[r] & \underset{}{\overset{-2}{\circ}}  \\
  \underset{}{\overset{-1}{ \bullet}} \ar@{-}[r] &    \underset{}{\overset{-2}{\circ}}  \ar@{-}[r]   & \underset{}{\overset{-3}{\circ}}  
&                  & \\
}
\]
\item If $H$ is of type $D_5$, then
\[
(Q\in H_Y) \cong (x^2+y^4+z^6-y^2z^2=0)/\mathbb{Z}_2(1,1,1)
\]
and the extended dual graph of $C \subset H_Y$ is
\[
\xymatrix @R=-7pt @C=5pt{
\underset{}{\overset{-1}{\bullet}} \ar@{-}[r]   &\underset{}{\overset{-2}{\circ}} \ar@{-}[dr] &
&                    & \underset{}{\overset{-2}{\circ}} \ar@{-}[dl]\\
                 &   &\underset{}{\overset{-3}{\circ}} \ar@{-}[r]\ar@{-}[dl]          & \underset{}{\overset{-2}{\circ}} \ar@{-}[dr]   & \\
                 &   \underset{}{\overset{-2}{\circ}}                                 &                                   &
& \underset{}{\overset{-2}{\circ}}
}
\]
\item If $H$ is of type $D_4$, then
\[
(Q \in H ) \cong (x^2+y^4+yt^3=0)/\mathbb{Z}_2(1,1,1)
\]
and the extended dual graph of $C \subset H_Y$ is
\[
\xymatrix @R=-1pt @C=5pt{
 & \overset{-2}{\circ} \ar@{-}[d] & &\\
 \underset{}{\overset{-2}{\circ}} \ar@{-}[r] & \circ \ar@{-}[r] \ar@{-}[d] & \underset{}{\overset{-2}{\circ}} \ar@{-}[r]
& \underset{}{\overset{-1}{\bullet}}\\
 & \underset{-2}{\circ} & &
}
\]
wwhere the unmarked curve has self intersection $-3$.
\end{enumerate}
\end{enumerate}
\end{theorem}
The proof of Theorem~\ref{E7} is given in subsection 4.1.

\begin{example}
Let $X$ be given by $x^2+y^3+yz^3+t^n=0$, and $\Gamma$ the line $x=y=t=0$. The general section through $\Gamma$ is an $E_7$ singularity and by
Proposition~\ref{normal-forms},
the general section of $X$ is $E_7$ if $n\geq 5$, $E_6$ if $n=4$ and $D_4$ if $n=3$. Theorem~\ref{E7} shows that a terminal
divisorial extremal neighborhood $C \subset Y \stackrel{f}{\la} X \supset P$ contracting an irreducible divisor $E$ onto $\Gamma$ exists
if and only if $n=3$.
\end{example}

In Theorem~\ref{E7}, the conditions for the existence of a terminal divisorial extremal neighborhood 
are given in terms of the form of the equations of $P\in \Gamma \subset X$.
The next proposition gives another criterion that is independent of equations.
\begin{proposition}\label{log-criterion}
Let $P \in \Gamma \subset X$ be the germ of a Gorenstein terminal 3-fold singularity $X$ along a smooth curve $\Gamma$. Suppose that the general
hyperplane section $S$ of $X$ through $\Gamma$ is DuVal of type $E_7$. Let $T \in |-2K_X|$ be the general member through $\Gamma$. Then
a terminal divisorial extremal neighborhood $C \subset Y \stackrel{f}{\la} X \supset P$ contracting an irreducible divisor $E$ onto $\Gamma$ 
exists iff the pair $(X,1/2T)$ is canonical with $e(X,1/2T)=1$.
\end{proposition}

\begin{proof}
Suppose that a terminal divisorial extremal neighborhood $C \subset Y \stackrel{f}{\la} X \supset P$ contracting an irreducible 
divisor $E$ onto $\Gamma$ exists. Then by Theorem~\ref{E7} $Y$ has index 2. Hence $-2K_Y$ is
very ample. Let $T_Y \in |-2K_Y|$ be the general element. Then $T_Y$ is smooth and it avoids the singular point of $Y$. Hence $(Y,1/2T_Y)$ is terminal.
Let $T=f_{\ast}T_Y\in |-2K_X|$. Let $l$ be a general fiber of $f$ over $\Gamma$. Then $l \cdot T_Y =-2$ and hence generically over $\Gamma$, $T_Y \la T$
is $2$ to $1$. Hence $T$ has a double point at the generic point of $\Gamma$. Since $f$ is generically the blow up of $\Gamma$, it follows that \[
K_Y+1/2T_Y=f^{\ast}(K_X+1/2T)
\]
and therefore $(X,1/2T)$ is canonical with $e(X,1/2T)=1$ as claimed.

Conversely, assume that $(X,1/2T)$ is canonical with $e(X,1/2T)=1$. Let $E$ be the unique crepant divisor for the pair.
Then by~\cite{Ko92} there exists an
extraction of $E$, i.e., a contraction $f \colon E \subset Y \la X \supset \Gamma$ such that $K_Y+1/2T_Y=f^{\ast}(K_X+1/2T)$.
Hence $(Y,1/2T_Y)$ is terminal
and hence so is $Y$.
\end{proof}

\begin{remark}
The condition $e(X,1/2T)=1$ is essential as shown by the next example.
\end{remark}
\begin{example}
Let $X$ be given by $x^2+y^3+yz^3+yt^2+t^4=0$ and $\Gamma$ by $x=y=t=0$.
The section $H$ given by $z=0$ is a $D_4$ DuVal
singularity and it follows from Theorem~\ref{E7}.2.c that there is no terminal divisorial extremal neighborhood over $P\in \Gamma \subset X$. 
The general $T\in |-2K_X|$ is given by $y=t^2$ and it is
$x^2+t^6+t^4+t^2z^3+t^4=0$. Then $C= H \cap T$ is given by $(x^2+t^4+t^6=0)\subset \mathbb{C}^3$. Straightforward calculations show that $(H,1/2C)$
is log-canonical. Therefore by inversion of adjunction, $(X,H+1/2T)$ is log-canonical and hence since $H$ is Cartier, $(X,1/2T)$ is canonical.
\end{example}

Theorem~\ref{E7} suggests a relation between divisorial extremal neighborhoods $C \subset Y \stackrel{f}{\la} X \ni P $ and
deformations of the general member of $|\sheaf_Y|$. First we introduce some notation that will simplify the statements.
\begin{definition}
\begin{enumerate}
\item We denote by $C\subset H_{7,k}$, $k=4,\; 5,\; 6$, the surface germs $C \subset H_Y$ that are described in Theorem~\ref{E7}.2.a,~\ref{E7}.2.b 
and~\ref{E7}.2.c, respectively.
\item Let $C \subset Y \stackrel{f}{\la} X \ni P $ be a divisorial extremal neighborhood. 
Let $E$ be the exceptional divisor and $\Gamma = f(E)$ its center.
Suppose that $\Gamma$ is smooth and that the general member $S\in |-K_X|$ containing $\Gamma$ is DuVal of type $E_7$.
Let $H\in |\sheaf_X|$ be the general hyperplane section. Then:
\begin{enumerate}
\item The neighborhood is called of type $E_{7,6}$ if $H$ is DuVal of type $E_6$.
\item The neighborhood is called of type $E_{7,5}$ if $H$ is DuVal of type $D_5$.
\item The neighborhood is called of type $E_{7,4}$ if $H$ is DuVal of type $D_4$.
\end{enumerate}
\end{enumerate}
\end{definition}

\begin{proposition}\label{deformations}
Divisorial extremal neighborhoods of type $E_{7,k}$, $k=4,\; 5,\; 6$ are in one to one correspondence with terminal
$\mathbb{Q}$-Gorenstein smoothings $Y$ of the surface germs $C\subset H_{7,k}$
such that the general $S_Y\in |-K_Y|$ is DuVal of type $E_6$, $C\cong \mathbb{P}^1$.
\end{proposition}

Let $f\colon C\subset Y \la X\supset P$ be a divisorial extremal neighborhood. Let $E$ be the $f$-exceptional divisor and $\Gamma$ its center in $X$. 
Assume that
the general hyperplane section $S$ of $X$ containing $\Gamma$ is DuVal of type $E_7$. Then Theorem~\ref{E7} and Proposition~\ref{deformations} show that
the general member of $|\sheaf_Y|$ is completely determined by the general member of $|\sheaf_X|$. It would be of interest to know if this holds in 
general. If it does then
one would be able to find the general member of any $|\sheaf_Y|$ by working out just a special example.
\begin{question}
Let $f\colon C\subset Y \la X\supset P$ be a divisorial extremal neighborhood. Let $E$ be the $f$-exceptional divisor and $\Gamma$ its center in $X$.
Is the general member of $|\sheaf_Y|$ completely determined by the general member of $|-K_X|$ and the general member of $|-K_X|$ 
containing $\Gamma$ (or equivalently
the general member of $|-K_Y|$)?
\end{question}
\begin{remark}
Surface germs $C\subset Z$ as in the previous proposition have terminal $\mathbb{Q}$-Gorenstein smoothings $Y$ that produce divisorial 
extremal neighborhoods
$f\colon C \subset Y \la X \ni P$ of type other that $E_{7,k}$, i.e., $Z$ is the general member of $|\sheaf_Y|$ but the general hyperplane 
section $S$ of $X$ containg
$\Gamma$ is not $E_7$. The type of the extremal neighborhoods associated to $\mathbb{Q}$-Gorenstein smoothings $Y$ of $Z$ is determined by 
the general member of $|-K_Y|$.
The way to get smoothings $Y$ with prescribed general anticanonical members is the following. Let $Q$ be the singular point of $Z$.
Let $\pi \colon Z^{\prime} \la Z$ be the index 1 cover of $Q\in Z$. This is given by an equation of the form $g(x,y,z)=0$ and it 
has a $\mathbb{Z}_2$-action. Let $Y^{\prime}$
be a smoothing given by $g(x,y,z)+tf(x,y,z,t)=0$ such that the $\mathbb{Z}_2$-action extends. In addition choose it in such a way 
such that the general hyperplane section
$S^{\prime}$ in the eigenspace of $-K_{Y^{\prime}}$ is DuVal. Then $S=S^{\prime}/\mathbb{Z}_2$ is the general member of $|-K_Y|$ and it is 
also DuVal. All finite covers of DuVal singularities
have been classified by M. Reid~\cite{Rei87} and looking at the list one knows what the $S^{\prime}$ must be in order to get the desired $S$. 
Next the deformation $Q\in Y$ lifts to a
deformation $Y$ of $Z$ and extending the contraction $f\colon Z\la H$ to $Y$ we get an extremal neighborhood $f \colon Y \la X$ with the 
desired invariants.

The next example demonstrates this.
\end{remark}
\begin{example}
Let $C\subset Z$ be the germ of a surface along a curve $C\cong \mathbb{P}^1$ such that the singular locus of $Z$ is a point $Q$ 
such that $(Q\in Z) \cong (x^2+y^4+yz^3=0)/\mathbb{Z}_2(1,1,1)$.
In this example we will construct a divisorial extremal neighborhood $f\colon C \subset Y \la X \supset P$ 
such that the general section of $|\sheaf_Y|$ is $Z$,
the general member of $|\sheaf_X|$ is Duval of type $D_4$ and the general member of $|-K_Y|$ is DuVal of type $E_7$ if $n=3$, or $D_4$ if $n=2$.
Hence if $\Gamma$ is the center of the $f$-exceptional divisor, then the general hyperplane section of $X$ containing $\Gamma$ is DuVal
of type $D_4$ if $n=2$ or of type $E_7$ if $n=3$.

Let $f\colon Z \la H$ be the contraction of $C$ to a DuVal singularity of type $D_4$.
Let $Q \in Y$ be the smoothing of $Q\in Z$ given by
\[
(x^2+y^4+yz^3+t^n=0)/\mathbb{Z}_2(1,1,1,0)
\]
where $n\in \{2,3\}$.
Let $\pi \colon Y^{\prime} \la Y$ be the index 1 cover of $Y$. This is given by $x^2+y^4+yz^3+t^n=0$. 
Let $S^{\prime}$ be the section given by $z=0$. This is in the same eigenspace as $-K_Y^{\prime}$
and hence $S=S^{\prime}/\mathbb{Z}_2(1,1,0)$ belongs to $|-K_Y|$. I claim that $S$ is DuVal of type $D_4$ if $n=2$, or $E_7$ if $n=3$.
Indeed. $S^{\prime}$ is given by $x^2+t^n+y^4=0$. This is an $A_3$ singularity if $n=2$ and $E_6$ if $n=3$.
Now setting $u=x^2$, $v=y^2$ and $xy=w$ it follows that $S$ is given by $w^2-uv=u+v^2+t^n=0$. Hence $S$ is given by $w^2-vt^n-v^3=0$ 
and therefore it is of type $D_4$ if $n=2$ and
$E_7$ if $n=3$. In any case, $Q\in Y$ is terminal~\cite{Rei87}.

Let \[
\text{Def}(Z) \la \text{Def}(Q\in Z)
\]
be the natural map. It is known that it is smooth and in particular surjective~\cite{Ko-Mo92}. Therefore the local deformation 
$Q\in Y$ of $Q\in Z$ lifts to a global deformation $Y$ of $Z$.
Moreover, since $K_Y$ is $S_2$, $P\in S$ lifts to a section $S\in |-K_Y|$. Let $f \colon Y \la X$ be the extension of $f\colon Z \la H$ to $Y$. 
Then this is a divisorial extremal neighborhood
such that $Z$ is the general member of $|\sheaf_Y|$. Let $E$ be the $f$-exceptional divisor and $\Gamma$ its center in $X$. 
Then by construction, if $n=2$ then
both the general hyperplane section of $X$ and the general through $\Gamma$ are DuVal of type $D_4$. On the other hand, if $n=3$, 
then the general hyperplane section of $X$ is $D_4$
and the general section containing $\Gamma$ is $E_7$.
\end{example}
\begin{proof}[Proof of Proposition~\ref{deformations}]
We will only show (2). The other cases are treated similarly.

Let $C \subset Y \stackrel{f}{\la} X \ni P $  be an extremal neighborhood of type $E_{7,5}$. Let $E$ be the $f$-exceptional divisor 
and $\Gamma$ its center in $X$.
Then according to Theorem~\ref{E7}, $C \cong \mathbb{P}^1$ and  the general $Z \in |\sheaf_Y|$ is a surface germ with the properties claimed.

Conversely, let $C \subset Z$ be a surface germ with properties as in the proposition. Let $Y$ be a terminal 
$\mathbb{Q}$-Gorenstein smoothing of $Z$ such that
the general $S_Y\in |-K_Y|$ is DuVal of type $E_7$.  Let $f \colon C \subset Z \la H \ni P$ be the contraction of $C$. 
Then $P\in H$ is DuVal of type $D_5$.
Moreover, there is a morphism~\cite[Proposition 11.4]{Ko-Mo92}
\[
\text{Def}(Z) \la \text{Def}(H)
\]
and hence $f$ extends to $Y$. Let $f\colon Y\la X$ be the resulting contraction. Since $X$ is Gorenstein, $f$ cannot be flipping and 
hence it is divisorial.
Let $E$ the exceptional divisor and $\Gamma = f(E)$ its center. It remains to show that $\Gamma$ is smooth and that $Z$ is 
the general member of $|\sheaf_Y |$.

Let $m_P(\Gamma)$ be the multiplicity of $\Gamma$ at $P$. Then~\cite{Tzi05}
\[
m_P(\Gamma) = -\frac{(K_Z \cdot C)^2}{C^2}
\]
Straightforward calculations show that $K_Z\cdot C = -1/2$ and $C^2=-1/4$. Therefore $m_P(\Gamma)=1$ and hence $\Gamma$ is smooth.

It remains to show that $Z$ is the general member of $|\sheaf_Y|$. Let $S\in |-K_X|$ be the general member containg $\Gamma$. 
Then by adjunctions $K_{S_Y}=f^{\ast}K_S$.
Therefore $S$ is DuVal and it must be of either $E_7$ or $E_8$. It cannot be $E_8$ because $E_8$ is factorial and hence it does not 
contain any lines. Hence it is $E_7$.
Let $T\in |\sheaf_X|$ be the general hyperplane section of $X$. Then it must be either $D_5$, or $D_4$ or $A_n$, for some $n$. 
If it was $A_n$, then $S$ would also be of type $A_m$~\cite{Tzi03}.
Suppose that $T$ is of type $D_4$. Let $T_Y=f^{\ast}T$. Then by Theorem~\ref{E7}, $T_Y$ has exactly one singular point $Q$ and moreover
\[
(Q\in Z) \cong (x^2+y^4+z^6-y^2z^2=0)/\mathbb{Z}_2(1,1,1)
\]
Hence there is a $\mathbb{Q}$-Gorenstein deformation of $T_Y$ to $Z$. This deformation lifts to a deformation of the index 1 
covers~\cite{KoBa88}. The index 1 cover of $T_Y$ is
the cusp singularity given by $x^2+y^4+z^6-y^2z^2=0$. Its fundamental cycle is a cycle of two rational curves. 
The index 1 cover of $H$ is the simple elliptic
singularity given by $x^2+y^4+yt^3=0$. Then by~\cite[Theorem 9]{Karr77}, a deformation of $T_Y$ to $Z$ exists iff the Dynkin diagram 
of $T_Y$ is a subset of the Dynkin diagram
of $Z$. But this is not the case. Hence the general member of $|\sheaf_X|$ is a $D_4$ DuVal and hence $Z$ is the general member of $|\sheaf_Y|$.
\end{proof}

\subsection{Proof of Theorem~\ref{E7}}

According to Proposition~\ref{normal-forms}, in suitable local analytic coordinates the equations of $P \in \Gamma\subset X$ are
\[
x^2+y^3+yz^3+tf_2(y,t)+tf_{\geq 3}(y,z,t)=0
\]
with $f_{\geq 3}(0,z,0)=0$ and $\Gamma$ is given by $x=y=t=0$. In the affine chart $x=xt$, $y=yt$ the blow up $W$ of $X$ along $\Gamma$ is given by
\begin{equation}
x^2t+y^3t^2+yz^3+f_2(yt,t)+f_{\geq 3}(yt,z,t)=0
\end{equation}
Moreover, $E$ is given by $y=t=0$, $F$ by $z=t=0$ and hence $f^{-1}(\Gamma)=E+3F$ and the intersection $E\cap F$ is a line $L$ given by $y=z=t=0$. 
It is now clear that
in this chart $W$ has exactly one singularity on $L$. Similarly by checking the other two charts we see that $W$ has only one singular point, say $P$ 
on $L$. (There is also another
singular point $Q$ that lies in $F$ but not in $E$ but for our purposes it does not affect anything). Hence in order
to make $E$ $\mathbb{Q}$-Cartier, we only need to blow it up. Hence $Z$ is just the blow up of $W$ along $E$. Moreover, since $E\cap F$ is a line $L$ 
in $F$ and in particular Cartier,
$F_Z$ is isomorphic to $F$ since it is the blow up of $L$ in $F$. Hence $F_Z$ can be contracted to a point. Therefore there are no flips and $Y$ 
fits in the following diagram
\begin{equation}
\xymatrix{
            & Z \ar[dl]_{\nu} \ar[dr]^{\pi} &\\
W \ar[dr]_g &   & Y \ar[dl]^f \\
            & X & }
\end{equation}
Next we describe $Z$ explicitely. Let $C=\nu^{-1}(P)\cong \mathbb{P}^1$. In the affine chart $y=yt$, $Z$ is given by
\[
x^2+y^3t^4+yz^3+\frac{1}{t}f_2(yt^2,t)+\frac{1}{t}f_{\geq 3}(yt^2,z,t)=0
\]
Moreover, the exceptional curve $C$ is given by $x=z=t=0$ and as said earlier it is not contained in $F_Z$. The change of variables $y\mapsto y+a$ 
transforms the equation of $Z$ to
\[
x^2+(y+a)^3t^4+(y+a)z^3+\frac{1}{t}f_2((y+a)t^2,t)+\frac{1}{t}f_{\geq 3}((y+a)t^2,z,t)=0
\]
It is now clear that $Z$ is singular along $C$ unless $f_2(0,t)\neq 0$. Hence $Y$ is terminal if and only if $f_2(0,t)\neq 0$. Moreover from 
the previous
calculations it follows that $f^{-1}(0)_{red}=\mathbb{P}^1$. This together with
Proposition~\ref{normal-forms} shows the claims about when a terminal contraction exists.

Suppose now that $Y$ is terminal. Next we will find its index.

\textit{Claim:} $Y$ has index 2.

By adjunctions it follows that
\begin{gather*}
K_W=g^{\ast}K_X+E+3F\\
K_Z=\pi^{\ast}K_Y+aF_Z
\end{gather*}
for a positive number $a$. Let $l\subset F_Z=\mathbb{P}^2$ be a general line. Then
\[
l\cdot K_Z =l \cdot K_Y =-1
\]
Moreover and since $E_Z$ is Cartier, one can check by direct calculation that $E_Z \cdot l =1$. Moreover
\[
l\cdot (E_Z+3F_Z)=l\cdot (E+3F)= l\cdot K_Y =-1
\]
and hence
\[
F_Z\cdot l =-2/3
\]
Therefore $a=3/2$ and hence $Y$ has index 2. Moreover, $3F_Z$ is Cartier and one can check from the equations describing $Z$ that $3F_Z \cap C$ is a 
reduced point.
Hence $C \cdot F_Z=1/3$ and hence
\[
C \cdot K_Y =-1/2
\]

Let $S$ be the general section of $X$ through $\Gamma$. By assumption it is DuVal of type $E_7$. Its birational transform $S_W$ in $W$ is just the 
blow up of $S$ along $\Gamma$.
From the explicit description of $W$ it follows by a starightforward direct calculation that the intersection of $E$ and $S_W$ is a smooth curve 
that lies in the smooth part of $S_W$.
Therefore $S_Z$ is isomorphic to $S_W$ and in particular it does not contain $C$. This shows that $S_Y\cong S$, where $S_Y=f_{\ast}^{-1}S$.

It only remains to describe the general elements of $|\sheaf_Y|$. Let $H$ be the general section of $X$. This is given by $z=ax+by+ct$. Let $H_W$ 
be the birational transform of $H$ in $W$.
Let $D$ be the $g$-exceptional curve, i.e., $H_W\cap F$. In the affine chart $x=xt$, $y=yt$, $H_W$ is given by
\begin{gather*}
x^2+y^3t+yt^2(ax-by-c)^3+\frac{1}{t}f_2(yt,t)+\frac{1}{t}f_{\geq 3}(yt,t(ax-by-c),t)=0\\
zt(ax-by-c)
\end{gather*}
Moreover, from (5) it follows that $E$ is given by $y=t=0$. Since $Y$ is terminal, it follows from the previous discussion that $f_2(0,t)\neq 0$. 
Hence $H_W \cap E$
is given by $x^2=y=t=0$, and its support is a smooth point. Using $x$, $y$ as local analytic coordinates, $H_W=\mathbb{C}^2$,  $D$ is given by $x=0$ 
and $E \cap H_W$ is $x^2=y=0$.
Hence $H_Z$ is the blow up of an ideal $I_Q$ of $H_W$ supported at a smooth point $Q\in H_W$
and such that in local analytic coordinates $H_W=\mathbb{C}^2$ and $I_Q=(x^2,y)$. Let $U \la \mathbb{C}^2$ be the blow up of the ideal $(x^2,y)$ 
in $\mathbb{C}^2$. Let $E$ be the exceptional
divisor. Then straightforward calculations show that $E\cong \mathbb{P}^1$, $U$ has exactly one singular point $R$, $R\in U$ is DuVal of type $A_1$ 
and $E^2=-1/2$.

We will only do the case when $H$ is of type $E_6$. In some sense it is the most difficult because in this case $H_Y$ does not have log-canonical 
singularities. The other cases are treated
similarly.

So suppose that $H$ is of type $E_6$. Then $H_W$ is the blow up of $0\in H$. Let $C_W$ the $g$-exceptional curve. Then $H_W$ has exactly one 
singular point $P\in H_W$, $P\in H_W$ is of type $A_5$
and $C_W^2=-1/2$. Let $H_Z$ be the birational transform of $H_W$ in $Z$, $C_Z$ the birational transform of $C_W$ in $H_Z$ and $C$ the 
exceptional curve.
Then from our previous discussion, $H_Z$ has exactly two singular points. One $A_1$ type DuVal that lies on $C\cap C_Z$ and an $A_5$ 
on $C_Z$ but not on $C$. Moreover,
$C^2=-1/2$ and an easy calculation shows
also that $C_Z^2=-1$. Moreover, the extended dual graph of $C_Z \subset H_Z$ at $P$ is
\[
\xymatrix @R=5pt @C=5pt{
\circ \ar@{-}[r]   & \circ \ar@{-}[r] & \circ \ar@{-}[r] \ar@{-}[d] & \circ \ar@{-}[r] & \circ  \\
                   &                  & \bullet              &                  & \\
}
\]
Let $\phi \colon V \la H_Z$ be the minimal resolution. Then the extended dual graph is
\[
\xymatrix @R=5pt @C=5pt{
\circ \ar@{-}[r]   & \circ \ar@{-}[r] & \circ \ar@{-}[r] \ar@{-}[d] & \circ \ar@{-}[r] & \circ  \\
                   &    \circ  \ar@{-}[r]   & \bullet              &                  & \\
}
\]
where the marked one is the birational transform $C_V$ of $C_Z$ in $V$. A straighforward calculation shows that $C_V^2=-3$. 
All the other exceptional curves have self-intersection $-2$.

Now $\pi$ contracts $C_Z$ (because it lies in $F_Z$) and hence the extended dual graph of $C\subset H$ is
\[
\xymatrix @R=5pt @C=5pt{
\circ \ar@{-}[r]   & \circ \ar@{-}[r] & \circ \ar@{-}[r] \ar@{-}[d] & \circ \ar@{-}[r] & \circ  \\
  \underset{}{\overset{1}{ \bullet}} \ar@{-}[r] &    \circ  \ar@{-}[r]   & \underset{}{\overset{3}{\circ}}              &                  & \\
}
\]
as claimed.

It remains to describe $H_Y$ with equations. $H_Y$ has rational singularities and hence it is completely determined by the minimal resolution. Hence
it follows that $H_Y$ is the general member of any extremal neighborhood $f \colon C \subset Y \la X \ni 0$ such that the 
general member of $|-K_Y|$ is $E_7$
and the general section of $X$ is $E_6$. Hence we may work with an explicit example. So, let $\Gamma \subset X$ be given by
\[
x^2+(y+t)^3+yz^3+t^4=0
\]
where $\Gamma$ is given by $x=y=t=0$. The general hyperplane section $H$ of $X$ is an $E_6$ singularity and the general $S$ 
containing $\Gamma$ is $E_7$.
From our previous discussion, there exists a terminal contraction $f\colon E \subset Y \la X \supset \Gamma$ such that $Y$ 
has index 2. In this case there is an
alternative local construction of the contraction~\cite{Ko-Mo92}. Let $T_Y\in |-2K_Y|$ be the general element. This is smooth. 
Let $T=f(T_Y)\in|-2K_X|$.
Let $\pi \colon Y^{\prime} \la Y$ be the double cover corresponding to $T_Y$, and let $\nu \colon X^{\prime} \la X$ be the double 
cover corresponding to $T$.
Then there is a commutative diagram
\begin{equation}
\xymatrix{
Y^{\prime} \ar[r]^{\pi} \ar[d]_{f^{\prime}}   & Y\ar[d]^{f}\\
X^{\prime} \ar[r]^{\nu}        & X
}
\end{equation}
Moreover, $K_{Y^{\prime}}=\pi^{\ast}(K_Y+1/2T_Y)$, $K_{X^{\prime}}=\nu^{\ast}(K_X+1/2T)$ and as we have seen in the proof of Proposition~\ref{log-criterion},
$K_Y+1/2T_Y=f^{\ast}(K_X+1/2T_X)$. From these it follows that
\[
K_{Y^{\prime}}=\pi^{\ast}(K_Y+1/2T_Y)=\pi^{\ast}f^{\ast}(K_X+1/2T)={f^{\prime}}^{\ast}\nu^{\ast}(K_X+1/2T)={f^{\prime}}^{\ast}K_{X^{\prime}}
\]
and hence $f^{\prime}$ is crepant. Moreover, since $(Y,1/2T_Y)$ is terminal, $Y$ is also terminal~\cite[Proposition 5.20]{Ko-Mo98}. 
Hence $f^{\prime}$ is a crepant
divisorial contraction from a terminal to a canonical 3-fold. At this point, one could study such contractions and eventually 
classify the terminal ones.
We will only use this construction in order to describe the index 1 cover of $H_Y$ which is just $\pi^{\ast}H_Y$. As we have 
mentioned earlier, to do this it suffices to
consider a special example only. Hence let $X$ be the 3-fold given by
\[
x^2+(y+t)^3+yz^3+t^4=0
\]
and $\Gamma$ the curve given by $x=y=t=0$. According to Proposition~\ref{normal-forms}, the general hyperplane section $H$ of $X$ is 
of type $E_6$. So let $H$ be the
section given by $z=0$ and let $T$ be the section of $|-2K_X|$ given by $y=t^2$. This is given by
\[
x^2+t^2[(1+t)^3t+z^3+t^2]=0
\]
and one immediately sees that $T$ is nodal generically along $\Gamma$. With notation as in diagram (4.3), $X^{\prime}$ is the double cover of $X$ 
corresponding to $T$.
Hence
\begin{gather*}
X^{\prime}=\text{Spec} \frac{\mathbb{C}[x,y,z,t,s]}{(s^2-y+t^2,x^2+(y+t)^3+yz^3+t^4)}\cong \\
\text{Spec}\frac{\mathbb{C}[x,z,t,s]}{(x^2+(s^2+t^2+t)^3+(s^2+t^2)z^3+t^4)}
\end{gather*}
The $\mathbb{Z}_2$-action on $X^{\prime}$ is $(x,z,t,s)\mapsto (x,z,t,-s)$ and $X=X^{\prime}/\mathbb{Z}_2$.
Moreover, $\Gamma^{\prime}=\nu^{-1}(\Gamma)$ is given by $x=t=s^2=0$ and $X^{\prime}$ is an $A_1$ DuVal singularity generically at $\Gamma^{\prime}$.
A straightforward direct calculation shows also that the blow up of $X^{\prime}$ along $\Gamma^{\prime}$ is crepant and terminal and $Y^{\prime}$ 
is just the blow
up of $X^{\prime}$ along $\Gamma^{\prime}$. Hence $H_{Y^{\prime}}$ is just the blow up of $H^{\prime}$ at the origin, where $H^{\prime}=\nu^{\ast}H$.
We will now calculate all steps explicitly starting from $H$ in order to describe $H^{\prime}$.

$H$ is given by $z=0$. Hence $H^{\prime}$ is given by
\[
x^2+(s^2+t^2+t)^3+t^4=0
\]
In the chart $x=xs$, $t=ts$, $H_{Y^{\prime}}$ is given by
\[
x^2+(s+st^2+t)^3s+t^4s^2=0
\]
and the $\mathbb{Z}_2$-action is given by $(x,t,s)\mapsto (-x,-t,-s)$. The change of variables $s=us$, $x=ux$, where $u=1/(1+t^2)$, 
preserves the $\mathbb{Z}_2$-action and
the equation of $H_{Y^{\prime}}$ becomes
\[
x^2+(s+t)^3s+ut^4s^2=0
\]
After setting $s=vs$ and $t=vt$, where $v=1\sqrt{u}=\sqrt{1+t^2}$, the equation of  $H_{Y^{\prime}}$ becomes
\[
x^2+(s+t)^3s+t^4s^2=0
\]
All the change of variables preserve the $\mathbb{Z}_2$-action and hence
\[
H_Y \cong (x^2+(s+t)^3s+t^4s^2=0)/ \mathbb{Z}_2(1,1,1)
\]
as claimed.

\section{The $E_6$ case.}
\begin{theorem}\label{E6}
Let $P\in \Gamma \subset X$ be the germ of a 3-fold $X$ along a smooth curve $\Gamma$ such that the general hyperplane section $S$ of $X$ 
through $\Gamma$ is DuVal of type $E_6$.
Then the general hyperplane section $H$ of $X$ is either $D_4$, $D_5$ or $E_6$. Moreover,
\begin{enumerate}
\item Suppose $H$ is $D_4$. Then there is a terminal divisorial extremal neighborhood $ C \subset Y \stackrel{f}{\la} X \ni P$ contracting an 
irreducible divisor $E$ onto $\Gamma$.
\item Suppose that $H$ is $D_5$ or $E_6$. Then in suitable local analytic coordinates, $P\in X$ is given by
\[
x^2+xz^2+y^3+by^2t+tf_{\geq 3}(y,z,t)=0
\]
such that no monomial of the form $z^k$, $k\geq 4$, or $y^s$, $s \geq 3$, appears in $f_{\geq 3}(y,z,t)$.
Let $a_{i,j,k}$ be the coefficient of $y^iz^jt^k$ in $f_{\geq 3}(y,z,t)$. Then a terminal divisorial extremal neighborhood 
$ C \subset Y \stackrel{f}{\la} X \ni P$ contracting an irreducible divisor $E$ onto $\Gamma$ exists, if and only if
one of the following two conditions hold
\begin{enumerate}
\item \[
(a_{030}^2+a_{021})^2+a_{003}\neq 0 \]
\item \[
2a_{030}(a_{030}^2+a_{021})+a_{012}\neq 0\]
\end{enumerate}
\end{enumerate}
Suppose that a terminal divisorial extremal neighborhood $ C \subset Y \stackrel{f}{\la} X \ni P$ exists and let $S_Y \in |-K_Y|$ be the general member. 
Then $S_Y \cong S$, $Y$ has index $3$ and its singular locus consists of one terminal singularity of type $cD_4/3$
\end{theorem}

\begin{proof}
By using the normal forms obtained in Proposition~\ref{normal-forms} we will calculate explicitly all steps of the construction 
described in section 2.
The main difficulty is to describe $Z$ and the possible existence of flips. As was shown in section 2, $W$ is singular along $E\cap F$ 
and hence $Z$ is not simply the blow up of $E$ as in the $E_7$ case. However, we will show that $Z$ fits in the following diagram:
\begin{equation}\label{construction-of-Z}
\xymatrix{
    & W_2 \ar[dl]_{h_2} \ar[dr]^{h^{\prime}} & \\
 W_1 \ar[dr]_{h_1} &                                & Z \ar[dl]^h \\
                  & W                             & }
\end{equation}
where $W_1$ is the blow up of $W$ along $E$ and $W_2$ is the blow up of $W_1$ along the birational transform $E_1$ of $E$ in $W_1$. $h_1$ 
has one exceptional divisor $B_1$ and
$h_2$ is a small contraction. More precisely, $h_1^{-1}(E)=E_1+2B_1$ is Cartier, $E\cong E_1$ and $W_1$ has three singular points on $E_1$
if $H$ is $D_4$, two if $H$ is $D_5$ and one if $H$ is $E_6$. Hence the $\mathbb{Q}$-factorialization of $E_1$ is just the blow up $W_2$ 
of $W_1$ along $E_1$. Moreover, we will
show that if $E_2$, $F_2$ and $B_2$ are the birational transforms of $E$, $F$ and $B_1$ in $W_2$, then 
$E_2\cong E_1$, $F_2 \cong F \cong \mathbb{P}^2$ and no $h_2$-exceptional curve is contained in $B_2$. Hence
$Z$ is obtained from $W_2$ by just contracting $B_2$. Moreover, from the above construction follows that no $h$-exceptional curve is contained
in $F_Z$, and therefore $F_Z\cong F\cong \mathbb{P}^2$. Hence $F_Z$ can be contracted to a terminal singularity
and therefore no flips exist. Then in order to decide whether $Y$ is terminal or only canonical
we need to study the singularities of $Z$ away from $F_Z$ (or equivalently the singularities of $W_2$ along the $h_2$-exceptional curves and 
away from $B_2$, $F_2$.
If they are isolated terminal, then so is $Y$. If not, then
$Y$ is only canonical.

We now proceed to justify the above claims by calculating all steps explicitly by using normal forms. By Proposition~\ref{normal-forms}, in 
suitable local analytic coordinates,
$P\in X$ is given by
\[
x^2+xz^2+y^3+ayt^2+by^2t+tf_{\geq 3}(y,z,t)=0
\]
and $\Gamma$ by $x=y=t=0$. Let $f_3(y,z,t)=a_{030}z^3+a_{210}y^2z+a_{201}y^2t+a_{120}yz^2+a_{021}z^2t+a_{003}t^3+a_{102}yt^2+a_{012}zt^2+a_{111}yzt$
be the degree 3 part of $f_{\geq 3}(y,z,t)$.

We now calculate $W$, the blow up of $X$ along $\Gamma$. In the affine chart $x=xt$, $y=yt$, $W$ is given by
\begin{gather*}
x^2t+xz^2+y^3t^2+ayt^2+by^2t^2+a_{030}z^3+a_{210}y^2t^2z+a_{201}y^2t^3\\
+a_{120}yz^2t+a_{021}z^2t+a_{003}t^3+a_{102}yt^3+a_{012}zt^2+a_{111}yzt^2+f_{\geq 4}(yt,z,t)=0
\end{gather*}
Setting $t=0$ we see that $g^{-1}(\Gamma)=E+2F$, where $E$ is given by $x+a_{030}z=t=0$, and $F$ is given by $z=t=0$. Moreover, $W$ is 
singular along the curve $L=E\cap F$ given by
$x=z=t=0$. This is the main difference between the $E_6$ and the $E_7$ case, where $W$ had only one singular point.

Now let $h_1 \colon W_1 \la W$ be the blow up of $W$ along $E$. Over the previous chart, it is given by
\begin{gather*}
(x+a_{030}z)u-vt=0\\
x^2t+xz^2+y^3t^2+ayt^2+by^2t^2+a_{030}z^3+a_{210}y^2t^2z+a_{201}y^2t^3\\
+a_{120}yz^2t+a_{021}z^2t+a_{003}t^3+a_{102}yt^3+a_{012}zt^2+a_{111}yzt^2+f_{\geq 4}(yt,z,t)=0
\end{gather*}
in $\mathbb{C}^4_{x,y,z,t}\times \mathbb{P}^1_{u,v}$. In the affine chart $u=1$, $W_1$ is given by
\begin{gather*}
(vt-a_{030}z)^2+vz^2+y^3t+ayt+by^2t+a_{210}y^2tz+a_{201}y^2t^2\\
+a_{120}yz^2+a_{021}z^2+a_{003}t^2+a_{102}yt^2+a_{012}zt+a_{111}yzt+\frac{1}{t}f_{\geq 4}(yt,z,t)=0
\end{gather*}
Setting $t=0$ we see that $g_1^{-1}(E)=E_1+2B_1$, where $B_1$ is given by $z=t=0$ and $E_1$ by
\[
t=a_{030}^2+v+a_{120}y+a_{021}+zg(z)=0
\]
where
\[
\frac{1}{t}f_{\geq 4}(yt,z,t)=h(y,z,t)+z^3g(z)
\]
and $h(y,z,0)=0$.
A straightforward calculation shows that no $g_1$-exceptional curve is contained in $E_1$ and hence $E_1\cong E$. 
Moreover, from the discussion in section 3, at the generic point
of $L$, $W$ is DuVal given by $x^2+y^2z-z^4=0$ and $E$ is given by $x=z=0$. A straightforward calculation now shows that the blow 
up $W_1$ of $W$ along $E$ is smooth along $E_1$
over the generic point of $L$ and
hence the $\mathbb{Q}$-factorialization of $E_1$ is just the blow up $W_2$ of $W_1$ along $E_1$.

Next we describe $W_2$. Over the previous affine chart of $W_1$, it is given by
\begin{gather*}
st=w(a_{030}^2+v+a_{120}y+a_{021}+zg(z))\\
(vt-a_{030}z)^2+vz^2+y^3t+ayt+by^2t+a_{210}y^2tz+a_{201}y^2t^2\\
+a_{120}yz^2+a_{021}z^2+a_{003}t^2+a_{102}yt^2+a_{012}zt+a_{111}yzt+\frac{1}{t}f_{\geq 4}(yt,z,t)=0
\end{gather*}
in $\mathbb{C}^4_{v,y,z,t}\times \mathbb{P}^1_{s,w}$. In the affine chart $w=1$, $W_2$ is given by
\begin{gather*}
st=a_{030}^2+v+a_{120}y+a_{021}+zg(z)\\
v^2t-2a_{030}zv+sz^2+y^3+ay+by^2+a_{210}y^2z+a_{201}y^2t\\
+a_{003}t+a_{102}yt+a_{012}z+a_{111}yz+\frac{1}{t}h(y,z,t)=0
\end{gather*}
or equivalently by
\begin{gather}\label{W2}
(st-a_{030}^2-a_{120}y-a_{021}-zg(z))^2t-2a_{030}z(st-a_{030}^2-a_{120}y-a_{021}-zg(z))+sz^2\\ \notag
+y^3+ay+by^2+a_{210}y^2z+a_{201}y^2t
+a_{003}t+a_{102}yt+a_{012}z+a_{111}yz+\frac{1}{t}h(y,z,t)=0
\end{gather}
in $\mathbb{C}^4_{s,y,z,t}$. The $g_2$-exceptional set $C_2$ is over $B_1$ since $E_1$ is Cartier away from $B_1$. Therefore $C_2$ is given by
\[
z=t=y^3+ay+by^2=0
\]
and by Proposition~\ref{normal-forms} it follows that $C_2$ is a curve with three irreducible components if the general hyperplane 
section $H$ of $X$ is of type $D_4$, two
components if $H$ is of type $D_5$ and is irreducible if $H$ is of type $E_6$. Contract now the birational transform $B_2$ 
of $B_1$ in $W_2$ to get $Z$ as in
diagram~\eqref{construction-of-Z}. Moreover the above calculations show that $F_2 \cong F$ and therefore $F_Z \cong F \cong \mathbb{P}^2$.
Hence $F_Z$ can be contracted to a terminal singularity and there is a commutative diagram
\begin{equation}
\xymatrix{
            & Z \ar[dl]_{h} \ar[dr]^{\pi} &\\
W \ar[dr]_g &   & Y \ar[dl]^f \\
            & X & }
\end{equation}
From the discussion at the beginning of the proof it follows that $Y$ is terminal if and only if $W_2$ has isolated singularities on $C_2$. 
The linear term of~\eqref{W2} is
\[
[(a_{030}^2+a_{021})^2+a_{003}]t+[2a_{030}(a_{030}^2+a_{021})+a_{012}]z+ay
\]
and the claim about the existence of a terminal contraction follows immediately (one has to check the other affine charts too in order 
to see that the conditions obtained are necessary and sufficient as well).

It now remains to describe the general member of $|-K_Y|$ and find the index of the singularities of $Y$. 

In section 1 it was shown that at the generic point 
of $L=E\cap F$, $W$ is DuVal of type $D_5$, 
and in suitable local analytic coordinates, $W$ is given by $x^2+y^2z-z^4=0$, $E$ by $x=z=0$ and $F$ by $x-z^2=y=0$. 
Therefore $2E$ and $4F$ are Cartier at all but finitely many points. 
Hence $2E_Z$ and $4F_Z$ are Cartier. Now write
\[
K_Z=\pi^{\ast}K_Y+aF_Z
\]
for some positive $a \in \mathbb{Q}$. Let $l \subset F_Z\cong \mathbb{P}^2$ a general line. Then
\begin{equation}
a l \cdot F_Z=l \cdot K_Z =l \cdot K_W =-1
\end{equation}
Moreover, since $K_W=g^{\ast}K_X+E+2F$, and $K_Z=h^{\ast}K_W$, it follows that $(E_Z+2F_Z)\cdot l =-1 $ and hence 
\begin{equation}
F_Z \cdot l = -\frac{1}{2}(1+E_Z \cdot l)
\end{equation}
Now $E_Z\cdot l = 1/2 (2E_Z)\cdot l $. Since $2E_Z$ is Cartier, $2E_Z \cdot F_Z=mL$, where $L=E\cap F$ and $m \in \mathbb{Z}$. 
The exact value of $m$ can be computed 
at the generic point of $L$. More precisely
\[
m=\text{length} (2E_Z \cap F_Z)
\]
From the previous discussion it follows that
\[
2E_Z \cap F_Z = \frac {k[x,y,z]}{(y,x-z^2,x^2,z)}=k
\]
where $K=k(L)$ is the field of fractions of $L$. Therefore $m=1$ and hence $E_Z \cdot l = 1/2$. Now from  $(6.4)$ and $(6.5)$ it follows that 
$a=4/3$ and hence $K_Z=\pi^{\ast}K_Y+4/3 F_Z$. Since $4F_Z$ is Cartier, $Y$ has index $3$ as claimed. 

We now proceed to describe the general member $S_Y$ of $|-K_Y|$. Then $S=f_{\ast}S_Y$.
Let $S_W$ and $S_Z$ the birational transforms of $S$ in $W$, $Z$ and $Y$. Then there is a commutative diagram
\begin{equation}
\xymatrix{
            & S_Z \ar[dl]_{h} \ar[dr]^{\pi} &\\
S_W \ar[dr]_g &   & S_Y \ar[dl]^f \\
            & S & }
\end{equation}
A careful examination of the above diagram following the explicit description of the corresponding diagram for $Y$, shows that 
$h$ is an isomorphism and $\pi$ contracts the $g$-exceptional curve. Therefore, $f$ is an isomorphism and hence $S_Y \cong S$. 

Now let $\nu \colon \tilde{Y} \la Y$ be the index 1 cover of $Y$ around the high index point of $Y$. Since $Y$ has index $3$, $\nu$ 
is finite of degree $3$. Let $\tilde{S_Y}\in |-K_{\tilde{Y}}|$ be the general member. Then $ \tilde{S_Y} \la S_Y$ is 3-to-1 and since $S_Y$ 
is DuVal of type $E_6$, $\tilde{S_Y}$ is DuVal of type $D_4$~\cite{Rei87} and hence the high index point of $Y$ is a $cD_4/3$ type
compound DuVal singularity.

\end{proof}

\end{document}